\documentclass{elsart}

\usepackage{amssymb}
\newtheorem{theorem}{Theorem}

\newtheorem{remark}{Remark}
\newtheorem{proposition}{Proposition}

\newtheorem{definition}{Definition}
\begin{document}

\begin{frontmatter}





\thanks{Supported by  DFG 436 RUS
113/823/0-1 and   the Ministry of Education of the Russian
Federation, project 2.1.1/1399}

\title{Nonexistence results for  compressible non-Newtonian fluid
 with magnetic effects in the whole space}



\author{ Olga Rozanova}

\address{Department of Differential Equations \& Mechanics and Mathematics Faculty,
Moscow State University, Moscow, 119992,
 Russia}
\ead{rozanova@mech.math.msu.su}
\begin{abstract}We consider a generalization of the compressible barotropic Navier-Stokes
equations to the case of non-Newtonian  fluid in the whole space.
The viscosity tensor  is assumed to be coercive with an exponent
$q>1.$ We prove that if the total mass and momentum of the system
are conserved, then one can find a constant $q_0> 1$ depending on
the dimension of space $n$ and the heat ratio $\gamma$ such that for
$q\in [q_0,n)$ there exists no global in time smooth solution to the
Cauchy problem. We prove also an analogous result for solutions to
 equations of magnetohydrodynamic non-Newtonian  fluid in 3D space.

\end{abstract}

\begin{keyword}
compressible non-Newtonian fluid \sep classical solution \sep loss
of smoothness
\MSC 53Q30
\end{keyword}
\end{frontmatter}











In the present paper we study a system of equations for the velocity
$u(t,x): {\mathbb R}_+\times {\mathbb R}^n\to {\mathbb R}^n$ and
density $\rho(t,x): {\mathbb R}_+\times {\mathbb R}^n\to \overline
{\mathbb R}_+$ describing the motion of viscous compressible fluid
without taking into account of heat phenomena in the $n$ -
dimensional space. The system is written as
$$\partial_t \rho+{\rm div}_x (\rho u)=0,\eqno(1)$$
$$\partial_t(\rho u)+{\rm Div}_x (\rho u \otimes u)\,=\,{\rm Div}_x{\mathbb S},\eqno(2)$$
where $\mathbb S$ is the stress tensor. We denote ${\rm Div}$ and
$\rm div$ the divergency of tensor and vector, respectively. The
Stokes axioms of the movement of continuum \cite{Serrin} imply the
following representation of the stress tensor:
$$
{\mathbb S} =\sum\limits_{k=0}^{n-1} \,\alpha_k(\rho,J_1({\mathbb
D}),...,J_n({\mathbb D}))\,{\mathbb D}^k,\eqno(3)
$$
where
$${\mathbb D}=(D_{ij})=\frac{1}{2}\left( \frac{\partial u_i}{\partial x_j}+\frac{\partial u_j}{\partial x_i}\right)$$
is the shear rate tensor and $J_m({\mathbb D}),\,m=1,...,n,$ are its
invariants, $\alpha_k,\,k=0,...,n-1,\,$ are scalar functions.

In particular, for the classical model of Newtonian fluid the stress
tensor is given by
$$
{\mathbb S} =(\lambda\,{\rm div u}-p(\rho))\,{\mathbb I}+
2\mu\,{\mathbb D},\eqno(4)
$$
with constant coefficients of viscosity $\,\lambda,\,\mu\,$($\mu >
0, \,\lambda+\frac{2}{n}\mu >0$)
 and a
given nonnegative
function of pressure $p(\rho): {\mathbb R}_+\times
{\mathbb R}^n\to {\mathbb R}.$

In what follows we assume
$$
{\mathbb S} =-p(\rho){\mathbb I}+{\mathbb P}(\rho,\mathbb D),\quad
{\mathbb P}(\rho,\mathbb D)=\beta_0(\rho,{\rm div}u){\mathbb
I}\,+\,\beta(\rho,|{\mathbb D}|)\,{\mathbb D},\eqno(6)
$$
where $\,{\mathbb I}\,$ is the identity matrix, $\beta_0$ and
$\beta$ are smooth functions bounded at zero,  $|{\mathbb
D}|^2=\mathbb D:\mathbb D,$
$$
\mathbb A : \mathbb B =\sum_{i,j=1}^n\,{\mathbb A}_{ij}{\mathbb
B}_{ij}
$$
for all tensors $\mathbb A$ and $\mathbb B$ of rank two. The
pressure satisfies the state equation
$$p=A\rho^\gamma,\,\gamma={\rm const}>1,\,A={\rm const}>0.\eqno(7)$$
Moreover,  we suppose that $\mathbb P$ obeys the integral coercivity
condition, i.e.
$$\int\limits_{{\mathbb R}^n} \mathbb P(g,\mathbb B): \mathbb B \,dx\ge \nu\, \int\limits_{{\mathbb R}^n}|{\mathbb
B}|^{q} \,dx, \quad q={\rm const}>1,\quad \nu={\rm const}>0,\,
\eqno(8)$$ holds for all symmetric $\,(n\times n)\,$ matrices $\,
{\mathbb B}$ and all nonnegative  $g$.
 Evidently, condition
$$
\beta_0(g,{\rm tr}{\mathbb B})\,{\rm tr}{\mathbb B}\,+
\,\beta(g,|{\mathbb B}|)|{\mathbb B}|^2 \ge \nu\,|{\mathbb
B}|^{q}\eqno(8')
$$
implies (8).

In the present work we are going to find a condition on the exponent
$q,$ the heat ratio $\gamma$ and the dimension of space $n$ such
that  smooth solutions to system (1), (2), (6)--(8) considered as
laws of conservation  of mass and momentum in the whole space do not
exist globally in time.

The viscocity tensor ${\mathbb P}$ of form (6)  is the most widely
used in the models of motion of incompressible non-Newtonian fluids,
e.g. models due to Ostwald, de Waele, Spriggs, Eyring, Carreau,
Ellis, Bingham and many others (\cite{Necas} and references
therein). The model (6), (8') includes the power-law model as a
special subclass. Power-law models are by far the most commonly used
models for describing the response of a variety of fluids. Many
colloids and suspensions are very well described by such models
(\cite{Schowalter},\cite{Bird}), it is also a well accepted fact
that blood, in a homogenized sense, can be modeled as a power-law
fluid (see for example \cite{Cho}, \cite{Cokelet}). The power-law
models are quite popular in glaciology \cite{Hutter} and geology.

For compressible fluids most results concern with Newtonian fluid
(local in time existence of classical solutions, global existence of
weak solutions \cite{AntontsevKazhikhovMonakhov}, \cite{Lions}). The
modern state of art in the mathematical theory of non-Newtonian
fluid is described in \cite{MalekNecasRokytaNecas}, where
homogeneous incompressible fluids are studied in sufficient details,
while compressible fluids are hardly considered (only very weak
measure-valued solutions are obtained). For stronger classes, the
solvability of (1), (2) is proved for an arbitrary $n$ (for bounded
domain) in \cite{Mamontov1}, \cite{Mamontov2}, the paper
\cite{Mamontov3} is devoted  to the problem of improvement of the
smoothness of the solution obtained,  the global regularity
estimates were derived in dimensions one and two. Thus, there are
still many unsolved problems in the mathematical theory of
non-Newtonian fluids, and first of all this is the problem of being
globally well-posed.

Let us assume $\rho\ge 0$ and introduce the  total mass
$$
m=\int\limits_{{\mathbb R}^n}\rho \, d x,
$$
momentum
$$P=\int\limits_{{\mathbb R}^n}\rho u\, d x$$
and total energy
$$ \mathcal E(t)=\int\limits_{{\mathbb
R}^n}\left(\frac{1}{2}\rho
|u|^2+\frac{A\rho^\gamma}{\gamma-1}\right)\, d x
\,=E_{k}(t)+E_{i}(t),$$ where
 $\,E_k(t),\,$  and $\,E_i(t)\,$ are the
kinetic and internal components of energy, respectively. The
H\"older inequality implies that for a given solution the momentum
is finite provided the total mass and total energy are finite. To
obtain the finiteness of total mass and energy we have to specify
the   decay rate of  density and velocity as $\,|x|\to\infty.$

Namely, let us consider initial data
$$\rho_0(x)\in \mathbb L^1\cap \mathbb L^\gamma ({\mathbb R}^n),\quad \rho_0(x)\ge
0,\quad u_0(x)\in \mathbb H^q_1({\mathbb R}^n).\eqno(9)
$$

\begin{definition}
We say that a solution $(\rho,u)$ to the problem (1), (2), (7), (8),
(9) belongs to the class $\mathbb K_{q}$ if
$$\rho(t,x)\in C^1(0,T;\mathbb L^1\cap \mathbb
L^\gamma ({\mathbb R}^n)),\quad u(t,x)\in C^2(0,T;\mathbb
H^q_1({\mathbb R}^n)),\quad T>0.
$$
\end{definition}

\begin{proposition}If $\, n\ge 2\,$ and $\,q\in [q_0,n),\, q_0=\frac{2n\gamma}{n(\gamma-1)+2\gamma}$, then the total mass and
total energy are finite on the solutions of the class $\mathbb
K_{q}$.
\end{proposition}

{\it Proof.} It is evident that  the total mass and the internal
energy are finite. To prove the convergence of the integral of
kinetic energy we need the following inequality:
$$
\left(\int\limits_{{\mathbb R}^n}\,
|u|^{\frac{qn}{n-q}}\,dx\right)^{\frac{n-q}{n}} \le K
\,\int\limits_{{\mathbb R}^n}\, |Du|^q \,dx,\eqno(10)
$$
where the constant $K=\frac{q(n-1)}{2(n-q)}>0.\,$ The latter
inequality holds for $u\in H^q_1({\mathbb R}^n),\,$ $\,q\in (1,n),$
(\cite{Hebey}, p.22) and follows from the Sobolev embedding.
Further, let us note that if $\rho\in C^1(0,T;\mathbb L^1\cap
\mathbb L^\gamma(\mathbb R^n)),$ then $$\rho\in C^1(0,T;\mathbb
L^\sigma(\mathbb R^n)), \,\sigma\in (1, \gamma).$$ Indeed, the
H\"older inequality implies
$$\int\limits_{{\mathbb R}^n}\,\rho^\sigma \,dx\le
\left(\int\limits_{{\mathbb R}^n}\,\rho
\,dx\right)^{\frac{\gamma-\sigma}{\gamma-1}}\left(\int\limits_{{\mathbb
R}^n}\,\rho^\gamma \,dx\right)^{\frac{\sigma-1}{\gamma-1}}.\eqno(11)
$$
Thus, from (10) and (11) we obtain
$$2E_k(t)\le \|\rho\|_{\mathbb L^\sigma(\mathbb R^n)}\|u\|_{\mathbb L^{\frac{2\sigma}{\sigma-1}}(\mathbb R^n)}^2 \le
\|\rho\|_{\mathbb L^\sigma(\mathbb R^n)}\left(\int\limits_{{\mathbb
R}^n}\, |u|^{\frac{q_\sigma
n}{n-q_\sigma}}\,dx\right)^{\frac{\sigma-1}{\sigma}}\le$$
$${\rm const}\cdot \|\rho\|_{\mathbb L^\sigma(\mathbb R^n)}\,\left(\int\limits_{{\mathbb R}^n}\, |Du|^{q_\sigma}
\,dx\right)^
\frac{(\sigma-1)n}{(n-q_\sigma)\sigma},$$ where $\displaystyle
q_\sigma=\frac{2n\sigma}{(\sigma-1)n+2\sigma}.$ It is easy to check
that $q_\sigma$ decreases with $\sigma$ and $1<q_\sigma<n.$ This
completes the proof. $\square$

\begin{remark} As follows from Proposition 1, the total mass and energy
are finite for the solutions from the class $\mathbb K_q$, therefore
the total  momentum is finite as well. However, to guarantee a
natural property of conservation of total mass and momentum  we
should require additionally that $\rho, u, Du$ vanish as
$|x|\to\infty$  sufficiently fast at any fixed $t$. This allows to
eliminate the surface integral in the Stokes formula. In terms of
the Sobolev spaces it  signifies that we consider $\rho, u, Du$ from
the space $\mathbb H^2_l(\mathbb R^n),\,l>\frac{n}{2}+1$ (e.g.
\cite{Richtmyer}, Sec.5.13). Then  due to a sufficiently quick decay
of $\rho, u, Du$ as $\,|x|\to\infty \,$ and properties of the
viscosity tensor $\mathbb P$ the values of total mass $m$ and total
momentum $P$ are conserved. Nevertheless, more natural way is to
require the conservation of  $m$ and $P$ in advance. Then the fact
of conservation of these values dictates  a fast decay of variables
at infinity.
\end{remark}

\begin{theorem}
Assume $\, n\ge 2\,$ and $\,q\in (q_0,n)$. If the initial momentum
$P\ne 0$, then there exists no global in time  solution $(\rho,u)$
to the Cauchy problem (1), (2), (6)--(8), (9) from the class
$\mathbb K_q$ with conserved total mass and momentum.
\end{theorem}

{\it Proof.} We use property (8) to obtain the following estimate of
the total energy rate:
$$\mathcal E'(t)\le -\nu \int\limits_{{\mathbb
R}^n}\,|Du|^q\,dx.\eqno(12)
$$

The H\"older inequality implies
$$
|P|=\Big|\int\limits_{{\mathbb R}^n}\,\rho u \,dx\Big|\le
\left(\int\limits_{{\mathbb R}^n}\,\rho^{\frac{qn}{n(q-1)+q}}
\,dx\right)^{\frac{n(q-1)+q}{qn}}\left(\int\limits_{{\mathbb R}^n}\,
|u|^{\frac{qn}{n-q}} \,dx\right)^{\frac{n-q}{qn}}.\eqno(13)
$$
Further, using the Jensen inequality we have
$$
\left(\frac{1}{m}\,\int\limits_{{\mathbb
R}^n}\,\rho^{\frac{qn}{n(q-1)+q}}
\,dx\right)^\frac{(\gamma-1)(n(q-1)+q)}{n-q}\le
\frac{\int\limits_{{\mathbb R}^n}\,\rho^\gamma
\,dx}{m}=\frac{(\gamma-1) E_i(t)}{mA}\eqno(14)
$$
for
$$\frac{(\gamma-1)(n(q-1)+q)}{n-q}\ge 1.\eqno(15)$$

Thus, (13) and (14) give
$$
|P|\le K_1\, \left(E_i(t)\right)^{\frac{n-q}{qn(\gamma-1)}}
\left(\int\limits_{{\mathbb R}^n}\, |u|^{\frac{qn}{n-q}}
\,dx\right)^{\frac{n-q}{qn}},\eqno(16)
$$
with a positive constant $K_1$ that depends on $\gamma, n, m.$

Further, once more we use inequality (10). Namely, (12), (16), (10)
imply
$$
\mathcal E'(t)\le -\nu
\frac{|P|^q}{K_1^q\,K}\,(E_i(t))^{-\frac{n-q}{qn(\gamma-1)}}\le$$$$\le\,-\nu
\frac{|P|^q}{K_1^q\,K}\,(\mathcal
E(0))^{-\frac{n-q}{qn(\gamma-1)}}=const<0.
$$
The latter inequality contradicts to  nonnegativity of the total
energy of system. Inequality (16) can be rewritten as
$\,\gamma\ge\frac{qn}{n(q-1)+q}\,$ or $\,q>
q_1=\frac{n\gamma}{(n+1)(\gamma-1)+1}.$ One can see that $q_0>q_1$.
Thus, the proof is over. $\square$


Now we consider the flow of compressible non-Newtonian  magnetic
fluid in the space ${\mathbb R}^3.$ The governing system
 is a combination of the compressible
equations of non-Newtonian  fluid  and Maxwell's equations of
electromagnetism:
$$\partial_t \rho+{\rm div}_x (\rho u)=0,\eqno(17)$$
$$\partial_t(\rho u)+{\rm Div}_x (\rho u \otimes u)=({\rm curl}_x H)\times H+{\rm Div}_x{\mathbb S},\eqno(18)$$
$$\partial_t H-{\rm curl}_x (u\times H)=-{\rm curl}_x(\eta\,{\rm curl}_x H),\qquad {\rm div}_x H=0, \eqno(19)$$ where
the stress tensor $\mathbb S$ obeys the conditions (6), (7), (8) and
$H=(H_1,H_2,H_3)$ denotes the  magnetic field, $\eta$ is a
nonnegative constant. We restrict ourselves to the barotropic case.

Although the electric field $E$ does not appear in the MHD system
(17 -- 19), it is indeed induced according to the  relation
$E=\eta\,{\rm curl}_xH- u\times H$  by the moving conductive flow in
the magnetic field.

The results on the classical solvability of the magnetohydrodynamic
equations are very scanty even in the case of Newtonian flow. Let us
mention the local in time existence of solutions to the  Cauchy
problem for the density separated from zero  \cite{VolpertKhudiaev}.
What about global solvability, even for the one-dimensional case,
the global existence of classical solutions to the full perfect MHD
equations with large data remains unsolved (smooth global solutions
near the constant state in one-dimensional case are investigated in
\cite{KO}). The existence of global weak solutions in the Newtonian
case was proved recently \cite{HW2} (see also \cite{HW1} for
references).


Now we are going to prove an analog of Theorem 1 for non-Newtonian
MHD equations.


The total energy in this case is
$$\mathcal E=\int\limits_{\mathbb R^3}\left(\frac{1}{2}\rho
|u|^2+\frac{|H|^2}{2}+\frac{A\rho^\gamma}{\gamma-1}\right)\, d x
\,=E_{k}(t)+E_m(t)+E_{i}(t).$$ Thus, there arises a new component
$E_m(t),$ the magnetic energy.

\begin{definition}
We define the class $\mathbb K_{q}^H$ of solutions $(\rho,u,H)$ to
system (17 -- 19),  (7), (8) as
$$\rho(t,x)\in C^1(0,T;\mathbb L^1\cap \mathbb L^\gamma ({\mathbb R}^3)),
\quad u(t,x)\in C^2(0,T;\mathbb H^q_1({\mathbb R}^3)),$$
$$H(t,x)\in C^2(0,T;\mathbb L^2({\mathbb R}^3)),\quad T>0.
$$
\end{definition}

\begin{theorem}
Let  $\,q\in [\frac{6\gamma}{5\gamma-3},3).$ If the initial momentum
$P$ does not vanish, then there exists no global in time  solution
$(\rho,u, H)$ to the Cauchy problem (17)-- (19), (6) -- (8) with
data
$$\rho_0(x)\in \mathbb L^1\cap \mathbb L^\gamma ({\mathbb R}^3),\quad u_0(x)\in \mathbb H^q_1({\mathbb
R}^3),\quad H_0(x)\in \mathbb L^2({\mathbb R}^3),
$$
from the class $\mathbb K_{q}^H$ with conserved total mass and
momentum.
\end{theorem}

{\it Proof.} As follows from Proposition 1, for  $\,q\in
[\frac{6\gamma}{5\gamma-3},3)$ the kinetic energy is finite for the
solution from $\mathbb K_{q}^H$, therefore total energy, total mass
and total momentum $P$ are finite as well. Moreover, the momentum
$P$ and mass $m$ are assumed to be  conserved.

 Further, we can estimate the derivative of the
total energy $\mathcal E$ as
$$\mathcal E'(t)\le -\eta \int \limits_{{\mathbb R}^3}|{\rm curl}_x H|^2\, dx-
\nu \int \limits_{{\mathbb R}^3}|D u|^q\, dx\le 0,\eqno(21)$$
therefore $\mathcal E(t)\le \mathcal E(0).$  We can rewrite
literally the proof of Theorem 1 to get the inequality $\mathcal
E'(t)\le - \,const <0,$ which implies the contradiction with
nonnegativity of total energy. $\square$

\begin{remark} In fact, Theorems 1 and 2 are  extensions
of results of \cite{RozJDE} and \cite{RozHYP} proved for the
Navier-Stokes ($n\ge 3$) and MGD ($n=3$) equations. For the
viscoelastic fluid it is possible to consider also $n=2,$ whereas
for the Newtonian fluid the problem of global in time existence of
smooth solutions with finite mass and energy is open in the
two-dimensional space. For  smooth initial data with a compact
support the nonexistence of global in time classical solutions in
the case of the Newtonian fluid follows from \cite{Xin} for any $n$.
\end{remark}


\begin{thebibliography}{99}



\bibitem{Serrin} J.Serrin, {\it Mathematical Principles of Classical
Fluid Mechanics,} \,Handbuch f\"ur Physik, \,Bd.VIII, Berlin (1959).

\bibitem{Necas} $\rm J. M\acute{a}lek, J.Ne\check{c}as,
K.P.Rajagopal,\,$ Global analysis of the flows of fluids with
pressure-dependent viscocities, Arch.Rational Mech. Anal. 165 (2002)
243--269.

\bibitem{Schowalter} W.R.Schowalter, \,{\it Mechanics of non-Newtonian
Fluids.}\,Pergamon Press, 1978.

\bibitem{Bird}R.B.Bird, R.S.Amstrong, O. Hassager,\,{\it Dynamics of
Polymer Liquids},\, Jon Wiley \& Sons, New York,\,2nd edition, 1993.

\bibitem{Cho} Y.I.Cho, K.R.Kensey,\,Effects of the non-Newtonian
viscocity of blood on hemodynamics of diseased arterial flows.\,
Adv.Bioengineering 15 (1989), 147--148.

\bibitem{Cokelet}G.R.Cokelet,\, The rheology of human blood, \,In:{\it
Biomechanics:Its function and Objectives,} Editors Y.C.Fung,
N.Perrone \& M.Aliker, Prentice-Hall, 1972, 63 -- 103.

\bibitem{Hutter} K.Hutter,\,{\it Theoretical Glaciology}, Dordrecht:
D.Reidel, 1983.

\bibitem{AntontsevKazhikhovMonakhov}S.N.Antontsev, A.V.Kazhikhov,
V.N.Monakhov, \, Boundary problems of mechanics of inhomogeneous
liquids, \, Novosibirsk, 1983.


\bibitem{Lions} P.-L. Lions, \,Mathematical topics in fluid
mechanics, Vol.2, \, Clarendon Press, Oxdord, 1996.



\bibitem{MalekNecasRokytaNecas} J.Malek, J.Necas, M.Rokyta,
M.$\rm R\dot{u}\check{z}i\check{c}ka$,\,{\it Weak and measure-valued
solutions to evolutionary PDEs}\,Chapman \& Hall, London, 1996.

\bibitem{Mamontov1} A.E.Mamontov,\, Global solvability of
multidimensional Navier-Stokes equations of compressible nonlinearly
viscous fluid. I", Siberian Math.J, 40 (1999), 408 -- 420.

\bibitem{Mamontov2}A.E.Mamontov,\, Global solvability of
multidimensional Navier-Stokes equations of compressible nonlinearly
viscous fluid. II", Siberian Math.J, 40 (1999), 635 -- 649.


\bibitem{Mamontov3}A.E.Mamontov,\, Global regularity estimates for
multidimensional equations of compressible non-Newtonian fluids,
\,Mathematical Notes, 68 (2000), 312--325.

\bibitem{Hebey} E.Hebey,\, Sobolev spaces on Riemannian manyfolds.
Lecture Notes in Mathematics. Springer Berlin/ Heidelberg. Vol.1635,
1996.

\bibitem {Richtmyer}R.D. Richtmyer,  Principles of advanced mathematical physics, Vol.1, Springer  New
York/ Heiderberg/ Berlin, 1978.

\bibitem{VolpertKhudiaev} A.I.Volpert, S.I.Khudiaev,  \, On the
Cauchy problem for composite systems of nonlinear equations,\,
Mat.Sbornik \, 87 (1972), N4, 504 -- 528.

\bibitem {KO} S. Kawashima, M. Okada, Smooth global solutions for the
one-dimensional equations in magnetohydrodynamics, Proc. Japan Acad.
Ser. A Math. Sci. 58 (1982) 384 -– 387.

\bibitem{HW2} X. Hu, D. Wang, Global solutions to the three-dimensional full
compressible magnetohydrodynamic flows, Comm. Math. Phys. 283
(2008), 255 –- 284.

\bibitem{HW1} X. Hu, D. Wang, \, Compactness of weak solutions to the three-dimensional
compressible magnetohydrodynamic equations, J. Differential
Equations 245 (2008) 2176 -– 2198.


\bibitem{RozJDE} O. Rozanova,\,
Blow up of smooth solutions to the compressible Navier-Stokes
equations with the data highly decreasing at
infinity\,J.Diff.Equat., 245 (2008), 1762 -- 1774.


 \bibitem{RozHYP} O. Rozanova,\,Blow up of smooth solutions to
the barotropic compressible magnetohydrodynamic equations with
finite mass and energy, \,Proceedings of Symposia in Applied
Mathematics 2009. Hyperbolic Problems: Theory, Numerics and
Ap-plications. American Mathematical Society. Volume: 67 (2009), 911
-- 918.

\bibitem{Xin} Z.P.Xin, \, Blowup of smooth solutions to the compressible
Navier-Stokes equation with compact density,\, Comm.Pure Appl.Math.
{ 51}(1998), 229 -- 240.




\end{thebibliography}
\end{document}